\documentclass[twoside,10pt]{amsart} 
\usepackage{amsfonts,amssymb,latexsym,amsthm, epsfig,amsmath,amscd}
\newcommand{\mf}{\mathfrak}

\newcommand{\lra}{\longrightarrow}

\newcommand{\ol}{\overline}

\newcommand{\ra}{\rightarrow}

\newcommand{\mc}{\mathcal}
\newcommand{\mbb}{\mathbb}
\newcommand{\tn}{\textnormal}

\newcommand{\pf}{{\bf Proof : }}

\newtheorem{de}{Definition}[section]
\newtheorem{re}[de]{Remark}
\newtheorem{pr}[de]{Proposition} 
\newtheorem{tr}[de]{Theorem}
\newtheorem{lm}[de]{Lemma} 

\newtheorem{nt}[de]{Notation} 

\newtheorem{co}[de]{Corollary}

\newcommand{\LMD}{\Lambda}
\newcommand{\lmd}{\lambda}

\def\vp{\rm \vspace{0.2cm}}
\def\M{\rm Max}

\def\hb{\hfill$\Box$}

\def\GL{\rm GL}

\def\GQ{\rm GQ}
\def\SL{\rm SL}
\def\EQ{\rm EQ}

\def\E{\rm E}

\def\Ma{\rm M}
\def\G{\rm G}
\def\GQ{\rm GQ}

\def\Sp{\rm Sp}

\def\k{\rm K_1}
\def\K{\rm K}
\def\KH{\rm KH}

\def\GQ{\rm GQ}
\def\u{\rm U}
\def\I{\rm I}

\def\O{\rm O}

\def\C{\rm C}

\def\w{\rm W}
\def\sk{\rm SK_1}
\def\nk{\rm NK_1}
\begin{document}
 
\title{${\bf {\rm NK}_1}$ of Bak's unitary group over Graded Rings}

\author{Rabeya Basu}
\address{Indian Institute of
Science Education and Research (IISER) Pune,  India} 
\email{rabeya.basu@gmail.com, rbasu@iiserpune.ac.in}
\thanks{Research by the first author was supported by SERB-MATRICS grant for the financial
year 2020--2021. And, research by the second author was supported by IISER (Pune) post-doctoral research grant.  }
\author{Kuntal Chakraborty}
\address{Indian Institute of
Science Education and Research (IISER) Pune,  India} 
\email{kuntal.math@gmail.com}
\thanks{Corresponding Author: rabeya.basu@gmail.com, kuntal.math@gmail.com}

\date{}
\maketitle

{\small{\bf Abstract:}
For an associative ring $R$ with identity, 
we study the absence of $k$-torsion 
in ${\rm NK_1GQ}(R)$; Bass nil-groups 
for the general quadratic or Bak's unitary groups.
By using a graded version of Quillen--Suslin theory we deduce an analog for the graded rings. \vp \\



\section{\bf Introduction}

Let $R$ be an associative ring with identity element $1$. 
When $R$ is commutative, we define ${\sk}(R)$ as the kernel of 
the determinant map from the Whitehead group ${\k}(R)$ to the group of units of $R$. 
The Bass nil-group ${\nk}(R)=\tn{ker}
({\k}(R[X])\ra {\k}(R))$; $X=0$. {\it i.e.}, the subgroup consisting of elements $[\alpha(X)]\in {\k}(R[X])$ 
such that $[\alpha(0)]=[{\I}]$. Hence ${\k}(R[X])\cong {\rm NK}_1(R)\oplus {\k}(R)$. 
The aim of this paper is to study some properties of Bass nil-groups ${\rm NK_1}$ for the general quadratic groups or Bak's unitary groups.

It is well-known that for many rings, {\it e.g.} if $R$ is regular Noetherian, Dedekind domain, or any ring with finite global dimension, 
the group ${\rm NK}_1(R)$ is trivial. On the other hand, if it is non-trivial, then it is not finitely generated as a group. {\it e.g.} 
if $G$ is a non-trivial finite group, the group ring $\mbb{Z}G$ is not regular. In many such cases, 
it is difficult to compute ${\rm NK}_1(\mbb{Z}[G])$. In \cite{Harmon}, D.R. Harmon proved the triviality of this group when $G$ is finite group of square free order. 
C. Weibel, in \cite{weibel}, has shown the non-triviality of this group for $G$ = $\mbb{Z}/2\oplus\mbb{Z}/2$, $\mbb{Z}/4$, and $D_4$. 
Some more results are known for finite abelian groups from the work of R.D. Martin; {\it cf.}\cite{martin}. It is also known ({\it cf.}\cite{HL}) that for a general finite 
group $G$, ${\rm NK}_1(R[G])$ is a torsion group for the group ring $R[G]$. In fact, for trivial ${\rm NK}_1(R)$, every element of finite order of ${\rm NK}_1(R[G])$ 
is some power of the cardinality of $G$. For $R=\mbb{Z}$, this is a result of Weibel. In particular, if $G$ is a finite $p$-group ($p$ a prime), then every element
of  ${\rm NK}_1(\mbb{Z}[G])$ has $p$-primary order. In \cite{STI}, J. Stienstra showed that 
${\nk}(R)$ is a ${\w}(R)$-module, where ${\w}(R)$ is the ring of big Witt vectors ({\it cf.}\cite{BL1} and \cite{WEL1}).
Consequently, in (\cite{WEL}, \S 3), C. Weibel 
observed that if $k$ is a unit in $R$, then ${\sk}(R[X])$ has no 
$k$-torsion, when $R$ is a commutative local ring. Note that if $R$ is a commutative local ring then 
${\sk}(R[X])$ coincides with ${\nk}(R)$; indeed, if $R$ is a local ring then
${\SL}_n(R)={\E}_n(R)$ for all $n>0$. Therefore, we may replace $\alpha(X)$ by 
$\alpha(X)\alpha(0)^{-1}$ and assume that $[\alpha(0)]=[{\I}]$. In \cite{RB}, the first author 
extended Weibel's result for arbitrary associative rings. In this paper we prove the analog 
result for $\lambda$-unitary Bass nil-groups, {\it viz.} ${\rm NK_1GQ}^{\lambda}(R, \LMD)$, 
where $(R,\LMD)$ is the form ring as introduced by A. Bak in \cite{Bak1}. The main 
ingredient for our proof is an analog of Higman linearisation (for a subclass of 
Bak's unitary group) due to V. Kopeiko; {\it cf.}\cite{V}. 
For the general linear groups, Higman linearisation ({\it cf.}\cite{RB1}) allows us 
to show that ${\nk}(R)$ has a unipotent representative. The same result is not true in general 
for the unitary nil-groups. Kopeiko's results in \cite{V} explain a complete description 
of the elements of ${\rm NK_1GQ}^{\lambda}(R, \LMD)$ that have (unitary) unipotent representatives. 
Followings are the main results in this article. 
\begin{tr} \label{th1} 
\label{has no torsion}
	Let $[\alpha(X)]=\big[\begin{pmatrix}
	A(X)&B(X)\\C(X)&D(X)
	\end{pmatrix}\big]\in {\rm NK_1GQ}^{\lambda}(R,\Lambda)$ with \\
	$A(X)\in {\GL}_r(R[X])$ for some $r\in \mathbb{N}$. 
	Then $[\alpha(X)]$ has no $k$-torsion if $kR=R$.
\end{tr}

And, an analog for the graded rings: 
\begin{tr} \label{th3} 
	Let $R =R_0\oplus R_1\oplus \dots$ be a graded ring.
	Let $k$ be a unit in $R_0$. Let $N = N_0 +N_1 +\dots+N_r\in {\rm M}_r(R)$ be a nilpotent matrix, and ${\I}$ denote the identity matrix. If
	$[({\I}+N)]^k=[{\I}]$ in ${\rm K_1GQ}^{\lambda}(R, \Lambda)$, then $[{\I} +N] = [{\I} +N_0]$.
\end{tr}

In the proof of \ref{th3}, we have used a graded version of Quillen--Suslin's local-global principle 
for Bak's unitary group over graded rings. This unify and generalize the results proved in \cite{BRK}, \cite{RB}, \cite{BRK2}, and \cite{BS}.
\begin{tr}\label{LG principle} {\bf (Graded local-global principle)}
Let $R = R_0\oplus R_1 \oplus R_2 \oplus \cdots$ be a graded ring with identity $1$. 
	Let $\alpha\in {\GQ}(2n, R, \Lambda)$ be such that $\alpha\equiv {\I}_{2n}\pmod{R_+}$. If 
	$\alpha_{\mathfrak m}\in {\EQ}(2n, R_{\mathfrak m},\Lambda_{\mathfrak m})$, 
	for every maximal ideal $\mathfrak m\in \M(C(R_0))$, then $\alpha\in {\EQ}(2n,R, \Lambda).$
\end{tr}

\section{\bf Preliminaries} 

Let $R$ be an associative ring with identity element $1$. 
Let ${\Ma}(n, R)$ denote the additive group of $n\times n$ matrices, and 
${\GL}(n, R)$ denote the multiplicative group of $n\times n$ invertible matrices.  
Let $e_{ij}$ be the matrix with $1$ in the $ij$-th position and $0$'s
elsewhere. The elementary subgroup of ${\GL}(n, R)$ plays a key role 
in classical algebraic K-theory. We recall, 
\begin{de} \tn{{\bf Elementary Group ${\E}(n, R)$:} The subgroup of all 
matrices in ${\GL}(n, R)$ generated by $\{{\E}_{ij}(\lambda):\lambda \in R, i\ne j\}$, where 
${\E}_{ij}(\lambda)={\I}_n+\lambda e_{ij}$, and 
$e_{ij}$ is the matrix with $1$ in the $ij$-position and $0$'s
elsewhere.}
\end{de}

\begin{de} \tn{For $\alpha\in {\Ma}(r,R)$ and $\beta\in {\Ma}(s,R)$, the matrix  
$\alpha\perp \beta$ denotes its 
embedding in ${\Ma}(r+s, R)$ (here $r$ and $s$ are even integers in the non-linear cases), given by}
$$\alpha\perp \beta =
\left(\begin{array}{cc} \alpha & 0 \\ 0 & \beta
\end{array}\right).$$ 
\tn{There is an infinite counterpart: 
Identifying each matrix $\alpha\in {\GL}(n, R)$ 
with the large matrix $(\alpha\perp \{1\})$
gives an embedding of ${\GL}(n, R)$ into ${\GL}(n+1, R)$. 
Let ${\GL}(R)=\underset{n=1}
{\overset{\infty}\cup} {\GL}(n, R)$, and 
${\E}(R)=\underset{n=1}{\overset{\infty}\cup} {\E}(n, R)$ be the corresponding 
infinite linear groups.}
\end{de}
As a consequence of classical Whitehead Lemma ({\it cf.}\cite{Ba}) due to A. Suslin, one gets 
$$[{\GL}(R),{\GL}(R)]={\E}(R).$$ 
\begin{de} \tn{The quotient group} 
$${\k}(R)=\frac{{\GL}(R)}{[{\GL}(R),{\GL}(R)]}=\frac{{\GL}(R)}{{\E}(R)}$$
\tn{is called the {\bf Whitehead group} of the ring $R$. For }
\tn{$\alpha\in {\GL}(n, R)$, let $[\alpha]$ denote its equivalence class in 
${\k}(R)$.}
\end{de}

In the similar manner we define ${\k}$ group for the other types of classical groups; {\it viz.}, the symplectic Whitehead group 
${\k}{\Sp}(R)$ and the orthogonal Whitehead group ${\k}{\O}(R)$. 

This paper explores a uniform framework for classical type groups over graded structures. 
Let us begin by recalling the concept of {\it form rings} and {\it form parameter} as introduced by A. Bak in \cite{Bak1}. 
This allows us to give a uniform definition for classical type groups. 

\begin{de}
{\bf (Form rings):} {\tn Let $R$ be an associative ring with identity, and
with an involution $-:R\rightarrow R$, $a\mapsto \overline{a}$. Let $\lambda \in C(R)$ = the center of $R$,
with the property $\lambda\overline{\lambda}=1$. We define two additive subgroups of $R$, {\it viz.}}
$$\Lambda_{\rm max}=\{a\in R\mid a=-\lambda \overline{a}\}~ \textit{and}~ \Lambda_{\rm min}=\{a-\lambda\overline{a}\mid a\in R\}.$$
{\tn One checks that for any $x\in R$, 
$\Lambda_{\rm max}$ and $\Lambda_{\rm min}$ are closed under the conjugation operation
$a\mapsto \overline{x}ax$. }

\tn{A $\lambda$-form parameter on $R$ is an additive subgroup  
$\Lambda$ of $R$ such that $\Lambda_{\rm min}\subseteq \Lambda \subseteq \Lambda_{\rm max}$, 
and $\overline{x}\Lambda x\subseteq \Lambda$ for all $x\in R$. 
{\it i.e.}, a subgroup between two additive groups which is also closed  under the conjugation operation. 
A pair $(R,\Lambda)$ is called a form ring.}
\end{de}

To define Bak's unitary group or the general quadratic group, we fix a central element $\lambda \in R$ with $\lambda\ol{\lambda}=1$, and 
then consider the form $$\psi_n= \begin{pmatrix} 0 &  {\I}_n \\ \lmd {{\I}}_n &0\end{pmatrix}.$$ 
For more details, see \cite{RB}, and \cite{RB2}.  \vp

\noindent{\bf Bak's Unitary or General Quadratic Groups ${\GQ}$:} \vp 
 
$${\GQ}(2n, R,\LMD) ~ = ~ \{\sigma\in {\GL}(2n, R, \LMD)\,|\, \ol{\sigma}\psi_n
\sigma=\psi_n\}.$$

\subsection*{Elementary Quadratic Matrices :} Let $\rho$ be the permutation, defined by
$\rho(i)=n+i$ for $i=1,\dots,n$. For $a\in R$, and $1\leq i,j\leq n$, we define

\begin{center}
$q\varepsilon_{ij}(a)={\I}_{2n}+ae_{ij}-\overline{a}e_{\rho(j)\rho(i)}$ for $i\neq j$,

\[ qr_{ij}(a) = \left\{ \begin{array}{ll}
{\I}_{2n}+ ae_{i\rho(j)}-\lambda\overline{a}e_{j\rho(i)} & \text{for}~ i\neq j \\{\I}_{2n}+ae_{\rho(i)j} & \text{for}~ i=j
\end{array} \right. \]  \\
\[ ql_{ij}(a) = \left\{ \begin{array}{ll}
{\I}_{2n}+ ae_{\rho(i)j}-\overline{\lambda}\overline{a}e_{\rho(j)i} & \text{for}~ i\neq j \\{\I}_{2n}+ae_{\rho(i)j} & \text{for}~ i=j
\end{array} \right. \]  \\
\end{center}
(Note that for the second and third type of elementary matrices, if $i=j$, then
we get $a=-\lambda \overline{a}$, and hence it forces that $a\in\Lambda_{\rm max}(R)$. One checks that these
above matrices belong to $\GQ(2n,R,\Lambda)$; {\it cf.}\cite{Bak1}. \vp \\
\noindent{\bf $n$-th Elementary Quadratic Group} ${\EQ}(2n,R,\Lambda)$: \vp\\
The subgroup generated
by $q\varepsilon_{ij}(a),qr_{ij}(a) \text{and } ql_{ij}(a)$, for $a\in R$ and $1\leq i,j\leq n$. 
For uniformity we denote the elementary generators of ${\EQ}(2n,R,\Lambda)$ by $\eta_{ij}(*)$. \vp

\noindent{\bf Stabilization map:} There are standard embeddings:

\begin{center}
${\GQ}(2n, R, \LMD) \lra {\GQ}(2n+2, R, \LMD)$
\end{center} 
given by 

\begin{center}
$\begin{pmatrix}
  a & b \\
c & d 
 \end{pmatrix} \mapsto 
\begin{pmatrix}
a & 0 & b & 0 \\
0 & 1 & 0 & 0 \\
c & 0 & d & 0 \\
0 & 0 & 0 & 1 
\end{pmatrix}.$
\end{center}
Hence we have 
\begin{center}
 ${\GQ}(R,\LMD) = \underset{\lra}\lim\,\, {\GQ}(2n, R, \LMD)$.
\end{center}

It is clear that the stabilization map takes generators of ${\EQ}(2n,R,\Lambda)$ to the
generators of ${\EQ}(2(n + 1),R,\Lambda)$. Hence we have 
\begin{center}
  ${\EQ}(R,\LMD) = \underset{\lra}\lim \,\,{\EQ}(2n, R, \LMD)$
\end{center} \vp 

There are standard formulas for the commutators between quadratic 
elementary matrices. For details, we refer \cite{Bak1} (Lemma 3.16). 
In later sections there are repeated use of those relations.
The analogue of the Whitehead Lemma for the general quadratic groups ({\it cf.}\cite{Bak1}) due to Bak allows us to write: 
$$[{\GQ}(R,\LMD), {\GQ}(R,\LMD)]=[{\EQ}(R,\LMD), {\EQ}(R,\LMD)]={\EQ}(R,\LMD).$$
Hence we define the {\bf Whitehead group} of the general quadratic group 
$${\k}{\GQ}=\frac{{\GQ}(R,\LMD)}{{\EQ}(R,\LMD)}.$$
And, the Whitehead group at the level $m$ 
$${\K}_{1,m}{\GQ}  = \frac{{\GQ}_m(R,\LMD)}{{\EQ}_m(R,\LMD)},$$
where $m=2n$ in the non-linear cases. 

Let $(R,\Lambda)$ be a form ring. We extend the involution of $R$ to the ring $R[X]$ of polynomials by setting $\overline{X}=X$. 
As a result we obtain a form ring $(R[X],\Lambda[X])$.

\begin{de}
\tn{The kernel of the group homomorphism} 
$${\rm K_1GQ}(R[X],\Lambda[X])\rightarrow {\rm K_1GQ}(R,\Lambda)$$
\tn{induced from  the  form ring homomorphism $(R[X],\Lambda[X])\rightarrow (R,\Lambda): X\mapsto 0$ is denoted by 
${\rm NK_1GQ}(R,\Lambda)$. We often say it as Bass nilpotent unitary ${\k}$-group of $R$, or just unitary nil-group.}
\end{de}

From the definition it follows that $${\rm K_1GQ}(R[X],\Lambda[X])={\rm K_1GQ}(R,\Lambda)\oplus {\rm NK_1GQ}(R,\Lambda).$$

In this context, we will use following two types of localizations, mainly over graded ring $R =  R_0\oplus R_1 \oplus R_2 \oplus \cdots$. 
\begin{enumerate}
\item Principal localization: for a non-nilpotent,
non-zero divisor $s$ in $R_0$ with $\ol{s}=s$, we consider the multiplicative subgroup $S = \{1,s,s^2,\dots\}$, and denote localized form ring by $(R_s, \Lambda_s)$. 
\item Maximal localization: for a maximal ideal $\mf{m}\in \M(R_0)$, we take the multiplicative subgroup $S=R_0- \mf{m}$, and 
denote the localized form ring by $(R_{\mf{m}}, \Lambda_{\mf{m}})$.
\end{enumerate}
\textbf{Blanket assumption}: We always assume that $2n\geq 6$. \vp

Next, we recall the 
well-known ``{\bf Swan--Weibel's homotopy trick}'', which is the main ingredient to handle the graded case. Let 
$R = R_0\oplus R_1 \oplus R_2 \oplus \cdots$ be a graded ring. 
An element $a \in R$ will be denoted by $a = a_0 + a_1 + a_2 + \cdots $, where $a_i\in R_i$ for each $i$, 
and all but finitely many $a_i$ are zero. Let $R_+= R_1 \oplus R_2 \oplus \cdots$. 
Graded structure of $R$ induces a graded structure on ${\rm M}_n(R)$ (ring of $n \times n$ matrices).

\begin{de}
\tn{Let $a \in R_0$ be a fixed element. We fix an element $b = b_0 + b_1 + \cdots$ in $R$ and define a ring homomorphism 
$\epsilon: R \rightarrow R[X]$ given by \[ \epsilon(b) = \epsilon (b_0 + b_1 + \cdots )\; = \; b_0 + b_1X + b_2X^2 + \cdots + b_iX^i + \cdots .\]
Then we evaluate the polynomial $\epsilon(b)(X)$ at $X = a$ and denote the image by $b^+(a)$ {\it i.e.}, $b^+(a) = \epsilon(b)(a)$.
Note that $\big(b^+(x)\big)^+(y) = b^+(xy)$. Observe, $b_0=b^{+}(0)$. 
We shall use this fact frequently.}
\end{de}

The above ring homomorphism $\epsilon$ induces a group homomorphism at the ${\GL}(2n,R)$ level for every $n \geq 1$, {\it i.e.}, for 
$\alpha \in {\GL}(2n,R)$ we get a map $$\epsilon: {\GL}(2n,R,\Lambda) \rightarrow {\GL}(2n,R[X],\Lambda[X]) \text{ defined by} $$ 
$$\alpha = \alpha_0 \oplus \alpha_1 \oplus \alpha_2 \oplus \cdots \mapsto \alpha_0 \oplus \alpha_1X \oplus \alpha_2X^2 \cdots,$$
where $\alpha_i\in {\rm M}(2n,R_i)$. 
As above for $a \in R_0$, we define $\alpha^+(a)$ as $$\alpha^+(a) = \epsilon(\alpha)(a).$$

\begin{nt} \tn{
By ${\GQ}(2n,R[X],\Lambda[X],(X))$ we shall mean the group of all quadratic matrices over $R[X]$, 
which are ${\I}_{2n}$ modulo $(X)$. Also if $R$ is a graded ring, then by ${\GQ}(2n, R,\Lambda,(R_+))$ 
we shall mean the group of all quadratic matrices over $R$ which are ${\I}_{2n}$ modulo $R_+$.}
\end{nt}

The following lemma highlights very crucial fact which we use (repeatedly) in the proof of ``Dilation Lemma''.  

\begin{lm}
\label{injectivity_under_localisation}
Let $R$ be a Noetherian ring and $s\in R$. Then there exists a natural number $k$ such that the homomorphism 
${\GQ}(2n,R,\Lambda, s^kR)\rightarrow {\GQ}(2n,R_s,\Lambda_s)$ 
$($induced by localization homomorphism $R\rightarrow R_s)$ 
is injective. Moreover, it follows that the induced map ${\EQ}(2n,R,\Lambda, s^kR) \rightarrow {\EQ}(2n,R_s, \Lambda_s)$ is injective.
\end{lm}

For the proof of the above lemma we refer \cite{HV}, Lemma 5.1. Recall
that any module finite ring $R$ is direct limit of its finitely generated
subrings. Also, 
${\G}(R, \LMD) = \underset{\lra}\lim \, {\G}(R_i, \LMD_i)$, 
where the limit is taken over all finitely generated subring of $R$.
Thus, one may assume that $C(R)$ is Noetherian. 
Hence one may consider module finite (form) rings $(R, \LMD)$ with identity. \vp 

Now we recall few technical definitions and useful lemmas.

\begin{de}
\tn{A row $(a_1,a_2,\dots, a_n)\in R^n$ is said to be unimodular if there exists $(b_1,b_2,\dots,b_n)\in R^n$ 
	such that $\sum_{i=1}^{n}a_ib_i=1$. The set of all unimodular rows of length $n$ is denoted by ${\rm Um}_n(R)$. }
\end{de}

For any column vector $v\in (R^{2n})^t$ we define the row vector $\widetilde{v}=\overline{v}^t\psi_n$.

\begin{de}
\tn{We define the map $M: (R^{2n})^t\times (R^{2n})^t\rightarrow M(2n,R)$ and the inner product $\langle,\rangle$ as follows:}
\begin{align*}
M(v,w) &= v.\widetilde{w}-\overline{\lambda}\,\overline{w}.\widetilde{v}\\
\langle v,w\rangle &= \widetilde{v}.w
\end{align*}
\end{de}
\noindent Note that the elementary generators of the groups ${\EQ}(2n,R,\Lambda)$
are of the form ${\I}_{2n} +M(*_1,*_2)$ for suitably chosen standard basis vectors.

\begin{lm}
$(${\it cf.}\cite{Bak1}$)$ The group ${\E}(2n, R,\Lambda)$ is perfect for $n\geq 3$, {\it i.e.}, 
$$[{\EQ}(2n,R,\Lambda),{\EQ}(2n,R,\Lambda)]={\EQ}(2n,R,\Lambda).$$
\end{lm}
\begin{lm}
	\label{splitting lemma}
For all elementary generators of  ${\GQ}(2n,R,\Lambda)$ we have the following splitting property: for all $x,y\in R$,
$$\eta_{ij}(x+y)= \eta_{ij}(x)\eta_{ij}(y).$$ 
\end{lm}
${\pf}$ See pg. 43-44, Lemma 3.16, \cite{Bak1}.

\begin{lm}
	\label{product in group}
	Let $G$ be a group, and $a_i,b_i\in G$, for $i=1,2,\ldots,n$. 
	Then for $r_i=\Pi_{j=1}^{i}a_j$, we have $\Pi_{i=1}^{n}r_ib_ir_i^{-1}\Pi_{i=1}^{n}a_i$.
\end{lm}
\begin{lm}
	\label{quadratic plus elementary}
	The group ${\GQ}(2n,R,\Lambda,R_+)\cap {\EQ}(2n,R,\Lambda)$ generated by the 
	elements of the type $\varepsilon\eta_{ij}(*)\varepsilon^{-1}$, where $\varepsilon\in {\EQ}(2n,R,\Lambda)$ and $*\in R_+$.
\end{lm}
${\pf}$
Let $\alpha\in {\GQ}(2n,R,\Lambda,R_{+})\cap {\EQ}(2n,R,\Lambda)$. Then we can write 
$$\alpha= \Pi_{k=1}^{r}\eta_{i_kj_k}(a_k)$$ for some element $a_k\in R$, $k=1,\dots,r$. 
We can write $a_k$ as $a_k= (a_0)_k+(a_+)_k$ for some $(a_0)_k\in R_0$ and $(a_+)_k\in R_+$. 
Using Lemma \ref{splitting lemma}, we can write $\alpha$ as, $$\alpha = \Pi_{k=1}^{r}(\eta_{i_kj_k}(a_0)_k)(\eta_{i_kj_k}(a_+)_k).$$ 
Let $\epsilon_t= \Pi_{k=1}^{t}\eta_{i_kj_k}((a_0)_k)$ for $1\leq t\leq r$.
By the Lemma \ref{product in group}, we have 
$$\alpha = \left(\Pi_{k=1}^{r}\epsilon_k \eta_{i_kj_k}((a_+)_k)\epsilon_k^{-1}\right)\left(\Pi_{k=1}^{r}\eta_{i_kj_k}((a_0)_k)\right).$$

Let us write $A= \Pi_{k=1}^{r}\epsilon_k \eta_{i_kj_k}((a_+)_k)\epsilon_k^{-1}$ and $B=\Pi_{k=1}^{r}\eta_{i_kj_k}((a_0)_k)$. 
Hence $\alpha=AB$. Let `over-line' denotes the quotient ring modulo $R_+$. 
Now going modulo $R_+$, we have $\overline{\alpha}=\overline{AB}=\bar{A} \bar{B}=\overline{\I}_{2n}\bar{B}= \overline{\I}_{2n}$, 
the last equality holds as $\alpha\in {\GQ}(2n,R,\Lambda, R_+)$. Hence, $\overline{B}= \overline{{\I}}_{2n}$. 
Since the entries of $B$ are in $R_0$, it follows that $B={\I}_{2n}$. 
Therefore it follows that $$\alpha=\Pi_{k=1}^{r}\epsilon_k \eta_{i_kj_k}((a_+)_k)\epsilon_k^{-1}.$$ \hb
\section{\bf Quillen--Suslin Theory for Bak's Group over Graded Rings}
 
\subsection{Local--Global Principle}

\begin{lm}
	\label{key lemma}
	Let $(R,\Lambda)$ be a form ring and $v\in {\EQ}(2n,R,\Lambda)e_1$. Let $w\in R^{2n}$ be
	a column vector such that $\langle v,w\rangle=0$. Then ${\I}_{2n}+M(v,w)\in {\EQ}(2n, R,\Lambda)$.
\end{lm}

${\pf}$ Let $v=\varepsilon e_1$. Then we have ${\I}_{2n}+ M(v,w)= \varepsilon ({\I}_{2n}+M(e_1, w_1))\varepsilon^{-1}$, 
where $w_1=\varepsilon^{-1}w$. Since $\langle e_1,w_1\rangle = \langle v,w \rangle=0$, 
we have $w_1^T= (w_{11}, \dots, w_{1 n-1},0,\dots,w_{12n})$. Therefore, since $\lambda \bar{\lambda}=1$, 
we have $${\I}_{2n}+M(v,w)=\prod_{\substack{ 1\leq j\leq n \\ 1\leq i\leq n-1}}\varepsilon  ql_{in}(-\bar{\lambda}\overline{w}_{1 n+i}) q 
\varepsilon_{jn}(-\bar{\lambda} \overline{w}_{1j}) ql_{nn}^{-1}(*) \varepsilon^{-1}$$

\begin{lm}
	\label{alphaplusiselementary}
	Let $R$ be a graded ring. Let $\alpha\in {\EQ}(2n, R,\Lambda)$. Then for every $a\in R_0$ one gets $\alpha^+(a)\in {\EQ}(2n,R,\Lambda)$. 
\end{lm}
${\pf}$ Let $\alpha =\Pi_{k=1}^{t}({\I}_{2n}+aM(e_{i_k},e_{j_k})),$ 
where $a\in R$ and $t\geq 1$. Hence for $b\in R_0$, we have $\alpha^+(b)=\Pi_{k=1}^{t}({\I}_{2n}+a^+(b)M(e_{i_k},e_{j_k}))$. 
Now taking $v=e_i$ and $v=a^+(b)e_j$ we have $\langle v,w\rangle=0$ and ${\I}_{2n}+M(v,w)={\I}_{2n}+ 	a^+(b)M(e_i,e_j))$ 
which belongs to ${\EQ}(2n, R,\Lambda)$ by Corollary \ref{key lemma}. Hence we have $\alpha^+(b)\in {\EQ}(2n, R, \Lambda)$ for $b\in R_0$. \hb

\begin{lm}
\label{dilation}
{\bf (Graded Dilation Lemma)}
Let $\alpha\in {\GQ}(2n,R,\Lambda)$ with $\alpha^+(0)={\I}_{2n}$ and $\alpha_s\in {\EQ}(2n, R_s,\Lambda_s)$ 
for some non-zero-divisor $s\in R_0$. Then there exists $\beta\in {\EQ}(2n, R, \Lambda)$ 
such that $$\beta_s^+(b)=\alpha_s^+(b)$$ for some $b=s^l$ and  $l\gg 0$.
\end{lm}

${\pf}$ Since $\alpha_s\in {\EQ}(2n, R_s,\Lambda_s)$ with $(\alpha_0)_s={\I}_{2n}$, then $\alpha_s= \gamma$, 
where $\gamma_{ii}=1+g_{ii}$ where $g_{ii}\in (R^+)_s$ and $\gamma_{ij}=g_{ij}$ for $i\neq j$, 
where $g_{ij}\in (R^+)_s$. Choose $l$ large enough such that every denominator of $g_{ij}$ for all $i,j$ 
divides $s^l$. Then by Lemma \ref{alphaplusiselementary}, we have $\alpha_s^+(s^l)\in {\EQ}(2n, R_s,\Lambda_s)$. 
As all denominator is cleared then $\alpha_s^+(s^l)$ permits a natural pullback. 
Hence we have $\alpha^+(s^l)\in {\EQ}(2n,R,\Lambda).$ Call this pullback as $\beta$.\hb

\begin{lm}
	\label{dilation is true}
	 Let $\alpha_s\in {\EQ}(2n,R_s,\Lambda_s)$ with $\alpha_s^+(0)= {\I}_{2n}$. 
	 Then one gets $$\alpha_s^+(b+d)\alpha_s^+(d)^{-1}\in {\EQ}(2n, R, \Lambda)$$ 
	 for some $s,d\in R_0$ and $b=s^l, l\gg 0$. 
\end{lm}
${\pf}$ Since $\alpha_s\in {\EQ}(2n,R_s,\Lambda_s)$, we have $\alpha_s^+(X)\in {\EQ}(2n, R_s[X], \Lambda_s[X])$. 
Let $$\beta^+(X)= \alpha^+(X+d)\alpha^+(d)^{-1},$$ where $d\in R_0$. Then we have 
$$\beta^+_s(X)\in {\EQ}(2n, R_s[X], \Lambda_s[X])$$ 
and $\beta^+(0)= {\I}_{2n}$. Hence by Lemma \ref{dilation}, we have, there exists $b=s^l$, $l\gg 0$, 
such that $\beta^+(bX)\in {\EQ}(2n, R[X], \Lambda[X])$. Putting $X=1$, we get our desired result.\hb \vp

{\bf Proof of Theorem \ref{LG principle} -- Graded Local-Global Principle:} \vp \\
Since $\alpha_{\mathfrak{m}}\in {\EQ}(2n,R_{\mathfrak{m}},\Lambda_{\mathfrak{m}})$ for all $\mathfrak{m}\in {\rm Max}(C(R_0))$, 
for each $\mathfrak{m}$ there exists $s\in C(R_0)\setminus \mathfrak{m}$ such that $\alpha_s\in {\EQ}(2n,R_s,\Lambda_s)$. 
Using Noetherian property we can consider a finite cover of $C(R_0)$, say $s_1+\dots +s_r = 1$. From Lemma \ref{dilation}, 
we have $\alpha^+(b_i)\in {\EQ}(2n,R,\Lambda)$ for some $b_i=s_i^{l_i}$ with $b_1+\dots +b_r=1$. 
Now consider $\alpha_{s_1s_2\dots s_r}$, which is the image of $\alpha$ in $R_{s_1s_2\dots s_r}$. 
By Lemma \ref{injectivity_under_localisation}, $\alpha \mapsto \alpha_{s_1s_2\dots s_r}$ is injective. 
Hence we can perform our calculation in $R_{s_1s_2\dots s_r}$ and then pull it back to $R$. 
\begin{equation*}
\begin{split}
\alpha_{s_1s_2\dots s_r} = &\alpha_{s_1s_2\dots s_r}^{+}(b_1+b_2+\dots +b_r)\\
= &((\alpha_{s_1})_{s_2s_3\dots})^{+}(b_1+\dots +b_r)((\alpha_{s_1})_{s_2s_3\dots})^{+}(b_2+\dots +b_r)^{-1}\dots 
\\&((\alpha_{s_i})_{s_1\dots \hat{s_i}\dots s_r})^{+}(b_i+\dots +b_r)((\alpha_{s_i})_{s_1\dots \hat{s_i}\dots s_r})^{+}(b_{i+1}+\dots +b_r)^{-1}
\\&((\alpha_{s_r})_{s_1s_2\dots s_{r-1}})^{+}(b_r)((\alpha_{s_r})_{s_1s_2\dots s_{r-1}})^{+}(0)^{-1}
\end{split}
\end{equation*} 
Observe that $((\alpha_{s_i})_{s_1\dots \hat{s_i}\dots s_r})^+(b_i+\dots +b_r)((\alpha_{s_i})_{s_1\dots \hat{s_i}\dots s_r})^+
(b_{i+1}+\dots +b_r)^{-1}\in {\EQ}(2n,R,\Lambda)$ due to Lemma \ref{dilation is true} (here $\hat{s_i}$ 
means we omit $s_i$ in the product $s_1\dots \hat{s_i}\dots s_r$), and hence 
$\alpha_{s_1s_2\dots s_r}\in {\EQ}(2n, R_{s_1\dots s_r}, \Lambda_{s_1\dots s_r})$.  
This proves $\alpha\in {\EQ}(2n, R,\Lambda)$. \hb \vp 

\subsection{Normality and Local--Global}
Next we are going to show that if $K$ is a commutative ring with identity and
$R$ is an associative $K$-algebra such that $R$ is finite as a left $K$-module, then the
normality criterion of elementary subgroup is equivalent to the Local-Global
principle for quadratic group. (One can also consider $R$ as a
right $K$-algebra.)

\begin{lm}
	\label{semilocal}
	$(${\bf Bass;} {\it cf.}\cite{Ba1}$)$ Let $A$ be an associative $B$-algebra such that $A$ is finite
	as a left $B$-module and $B$ be a commutative local ring with identity. Then $A$ is semilocal.
\end{lm}

\begin{tr}
	\label{transitive action}
	$(${\it cf.}\cite{RB}$)$
	Let $A$ be a semilocal ring $($not necessarily commutative$)$ with involution. 
	Let $v\in {\rm Um}_{2n}(A)$.Then $v\in e_1{\EQ}(2n, A)$. In other words the group ${\EQ}(2n,A)$ acts transitively on ${\rm Um}_{2n}(A)$.
\end{tr}
Before proving the next theorem we need to recall a theorem from \cite{RB}:

\begin{tr}$(${\bf Local-Global Principle}$)$
	\label{lgfor quadratic} Let $A$ be an associative $B$-algebra such that $A$ is finite
	as a left $B$-module and $B$ be a commutative ring with identity..
	If $\alpha(X)\in {\GQ}(2n, A[X], \Lambda[X])$, $\alpha(0)=\rm {\rm I}_{2n}$ and $\alpha_{\mathfrak{m}}(X)\in 
	{\EQ}(2n, A_{\mathfrak{m}}[X],\Lambda_{\mathfrak{m}}[X])$ for every maximal ideal $\mathfrak{m}\in {\rm Max}(B)$, 
	then $\alpha\in {\EQ}(2n, A[X],\Lambda[X])$.
\end{tr}
\begin{tr} 
	Let $K$ be a commutative ring with unity and $R=\oplus_{i=0}^{\infty}R_i$ be a graded $K$-algebra such that $R_0$ is 
	finite as a left $K$-module. Then for $n\geq 3$ the following are equivalent:
	
	$(1)$ ${\EQ}(2n,R,\Lambda)$ is a normal subgroup of ${\GQ}(2n, R,\Lambda)$.
	
	$(2)$ If $\alpha\in {\GQ}(2n,R,\Lambda)$ with $\alpha^+(0)={\rm I}_{2n}$ and $\alpha_{\mathfrak{m}}\in {\EQ}(2n, R_{\mathfrak{m}}, 
	\Lambda_{\mathfrak{m}})$ for every maximal ideal $\mathfrak{m}\in {\rm Max}(K)$, then $\alpha\in {\EQ}(2n, R, \Lambda)$.
\end{tr}
	
${\pf}$ $(1)\Rightarrow (2)$ We have proved the Lemma \ref{key lemma} for any form ring with
identity and shown that the local-global principle is a consequence of Lemma \ref{key lemma}.
So, the result is true in particular if we have ${\EQ}(2n, R, \Lambda)$ is a normal subgroup of ${\GQ}(2n, R,\Lambda)$.

$(2)\Rightarrow (1)$ Since polynomial rings are special case of graded rings, the result
follows by using the Theorem \ref{lgfor quadratic}. Let $\alpha\in {\EQ}(2n, R,\Lambda)$ and $\beta \in {\GQ}(2n, R, \Lambda)$. 
Then we have $\alpha$ can be written as product of the matrices of the form  $({\rm I}_{2n}+\beta M(*_1,*_2)\beta^{-1})$, 
with $\langle*_1,*_2\rangle=0$ where $*_1$ and $*_2$ are suitably chosen basis vectors. Let $v=\beta *_1$. 
Then we can write $\beta \alpha \beta^{-1}$ as a product of the matrices of the form ${\rm I}_{2n}+M(v,w)$ 
for some $w\in R^{2n}$. We must show that each ${\rm I}_{2n}+M(v,w)\in {\EQ}(2n, R,\Lambda)$.

Consider $\gamma= {\rm I}_{2n}+M(v,w)$. Then $\gamma^+(0)={\rm I}_{2n}$. By Lemma \ref{semilocal} we have the ring $S^{-1}R$ 
is semilocal where $S=K\setminus \mathfrak{m}$, and $\mathfrak{m}\in {\rm Max}(K)$. Since $v\in {\rm Um}_{2n}(R)$, then by 
Theorem \ref{transitive action}, we have $v\in {\EQ}(2n, S^{-1}R,S^{-1}\Lambda)e_1$. 
Therefore by applying Lemma \ref{key lemma} to the ring $(S^{-1}R,S^{-1}\Lambda)$, 
we have $\gamma_{\mathfrak{m}}\in {\EQ}(2n, R_{\mathfrak{m}}, \Lambda_{\mathfrak{m}})$ 
for every maximal ideal $\mathfrak{m}\in {\rm Max}(K)$. Hence by hypothesis we have $\gamma\in {\EQ}(2n, R,\Lambda)$. This completes the proof. \hb

 \begin{re}
 	 \tn{We conclude that the local-global principle for the elementary subgroups and their normality properties are equivalent.}
 \end{re}

\section{\bf Bass Nil Group ${\rm {\nk}{\GQ}(R)}$} 

In this section recall some basic definitions and properties of the representatives 
of ${\rm {\nk}{\GQ}(R)}$.
We represent any element of
${\rm M}_{2n}(R)$ as {\small $\begin{pmatrix}
a&b\\c&d
\end{pmatrix},$} where $a,b,c,d\in {\rm M}_n(R)$. For {\small $\alpha =\begin{pmatrix}
a&b\\c&d
\end{pmatrix}$} we call {\small $\begin{pmatrix}
a&b
\end{pmatrix}$} the upper half of $\alpha$. Let $(R,\lambda, \Lambda)$ be a form ring. 
By setting $\bar{\Lambda}=\{\bar{a}: a\in \Lambda\}$ we get another form ring $(R,\bar{\lambda},\bar{\Lambda})$.
We can extend the involution of $R$ to ${\rm M}_n(R)$ by setting $(a_{ij})^*=(\overline{a}_{ji})$.

\begin{de}
	\tn{Let $(R,\lambda,\Lambda)$ be a form ring. A matrix $\alpha=(a_{ij})\in {\rm M}_n(R)$ 
	is said to be $\Lambda$-Hermitian if $\alpha=-\lambda \alpha^*$ and all the diagonal entries of $\alpha$ 
	are contained in  $\Lambda$. A matrix $\beta\in {\rm M}_n(R)$ is said to be $\bar{\Lambda}$-Hermitian 
	if $\beta=-\bar{\lambda}\beta^*$ and all the diagonal entries of $\beta$ are contained in $\bar{\Lambda}$. }
\end{de}

\begin{re}
	\tn{A matrix $\alpha\in {\rm M}_n(R)$ is $\Lambda$-Hermitian if and only if $\alpha^*$ is $\bar{\Lambda}$-Hermitian. 
	The set of all $\Lambda$-Hermitian matrices forms a group under matrix multiplication.  }
\end{re}

\begin{lm}
	\label{stableunderconjugation}
	\cite[Example 2]{V}
	 Let $\beta\in {\GL}_n(R)$ be a $\Lambda$-Hermitian matrix. Then the matrix $\alpha^*\beta \alpha$ is $\Lambda$-Hermitian for every $\alpha \in {\GL}_n(R)$.
\end{lm}

\begin{de}
	\label{lamdaunitary}
	\tn{Let $\alpha=\begin{pmatrix}
	a&b\\c&d
	\end{pmatrix}\in {\rm M}_{2n}(R)$ be a matrix. Then $\alpha$ is said to be a $\Lambda$-quadratic matrix if one of the following equivalent conditions holds:}
	\begin{enumerate} 
\item \tn{$\alpha \in {\GQ}(2n, R,\Lambda)$ and the diagonal entries of the matrices $a^*c, b^*d$ are in $\Lambda$,}
\item \tn{$a^*d+\lambda c^*d={\I}_n$ and the matrices $a^*c, b^*d$ are $\Lambda$-Hermitian,}
\item \tn{$\alpha\in {\GQ}(2n,R,\Lambda)$ and the diagonal entries of the matrices $ab^*,cd^*$ are in $\Lambda$,}
\item \tn{$ad^*+\lambda bc^*={\I}_{n}$ and the matrices $ab^*, cd^*$ are $\Lambda$-Hermitian. }
	\end{enumerate}
\end{de}

\begin{re}
\tn{The set of all $\Lambda$-quadratic matrices of order $2n$ forms a group called $\Lambda$-quadratic group. 
We denote this group by ${\GQ}^{\lambda}(2n, R, \Lambda)$. If $2\in R^*$, then we have $\Lambda_{\rm min}=\Lambda_{\rm max}$. 
In this case notions of quadratic groups and notions of $\Lambda$-quadratic groups coincides. Also this happens when $\Lambda=\Lambda_{\rm max}$.
Hence  quadratic groups are special cases of $\Lambda$-quadratic groups. Other classical groups appear as $\Lambda$-quadratic groups 
in the following way. Let $R$ be a commutative ring with trivial involution. Then 
\[
{\GQ}^{\lambda}(2n,R,\Lambda)=\begin{cases}
{\Sp}_{2n}(R), & \text{if } \lambda=-1 \text{ and } \Lambda=\Lambda_{\rm max}=R\\
{\O}_{2n}(R), & \text{if } \lambda=1 \text{ and } \Lambda=\Lambda_{\rm min}=0 
\end{cases}
\]
 And for general linear group ${\GL}_n(R)$, 
we have, ${\GL}_n(R)= {\GQ}^1(2n, H(R),\Lambda=\Lambda_{\rm max})$, where $\mbb{H}(R)$ denotes the ring 
$R\oplus R^{op}$ with $R^{op}$ is the opposite ring of $R$ and the involution on $\mbb{H}(R)$  
is defined by $\overline{(x,y)}=(y,x)$. Thus the study of $\Lambda$-quadratic matrices unifies the study of quadratic matrices.  }
\end{re}
 
We recall following results from \cite{V}. 
 
 \begin{lm}
 	Let $\alpha=\begin{pmatrix}
 	a&0\\0&d
 	\end{pmatrix}\in {\rm M}_{2n}(R)$. Then $\alpha\in {\GQ}^{\lambda}(2n, R,\Lambda)$ if and only if $a\in {\GL}_n(R)$ and $d=(a^*)^{-1}$.
 \end{lm}

${\pf}$ Let $\alpha \in {\GQ}^{\lambda}(2n, R, \Lambda)$. In view of $(2)$ of Definition \ref{lamdaunitary}, we have, $a^*d= {\I}_n$. Hence $a$ is invertible 
and $d=(a^*)^{-1}$. Converse holds by $(2)$ of Definition \ref{lamdaunitary}. \hb

\begin{de}
	\tn{Let $\alpha\in {\GL}_n(R)$ be a matrix. A matrix of the form {\small $\begin{pmatrix}
	\alpha&0\\0&(\alpha^*)^{-1}
	\end{pmatrix}$} is denoted by $\mbb{H}(\alpha)$ and is said to be hyperbolic.}
\end{de}

\begin{re}
	\tn{In a similar way we can show that matrices of the form  $T_{12}(\beta):=\begin{pmatrix}
	{\I}_n &\beta\\0&{\I}_n
	\end{pmatrix}$ is $\Lambda$-quadratic matrix if and only if $\beta$ is $\bar{\Lambda}$-Hermitian. And the matrix of the form $T_{21}(\gamma):=\begin{pmatrix}
	{\I}_n&0\\\gamma&{\I}_n
	\end{pmatrix}$ is $\Lambda$-quadratic matrix if and only if $\gamma$ is $\Lambda$-Hermitian.}
\end{re}

Likewise in the quadratic case we can define the notion of $\Lambda$-elementary quadratic groups in the following way:

\begin{de} \tn{
The $\Lambda$-elementary quadratic group is denoted by ${\EQ}^{\lambda}(2n, R, \Lambda)$ 
and defined by the group generated by $2n\times 2n$ matrices of the form $\mbb{H}(\alpha)$ 
where $\alpha\in {\E}_n(R)$, $T_{12}(\beta)$ and $\beta$ is $\bar{\Lambda}$-Hermitian and $T_{21}(\gamma)$ is $\gamma$ $\Lambda$-Hermitian.}
\end{de}

\begin{lm}
	\label{uppertriangular}
	Let $A=\begin{pmatrix}
	\alpha&\beta\\0&\delta
	\end{pmatrix}\in {\rm M}_{2n}(R)$. Then $A\in {\GQ}^{\lambda}(2n,R,\Lambda)$ if and only if 
	$\alpha\in {\GL}_n(R)$, $\delta=(\alpha^*)^{-1}$ and $\alpha^{-1}\beta$ is $\bar{\Lambda}$-Hermitian. In this case 
	$A\equiv \mbb{H}(\alpha) \pmod{{\EQ}^{\lambda}(2n, R, \Lambda)}$. 
\end{lm}

${\pf}$ Let $A\in {\GQ}^{\lambda}(2n, R, \Lambda)$. Then by $(4)$ of Definition \ref{lamdaunitary}, we have 
$\alpha\delta^*={\I}_n$ and $\alpha \beta^*$ is $\Lambda$-Hermitian. Hence $\alpha$ is invertible and 
$\delta=(\alpha^*)^{-1}$. For $\alpha^{-1}\beta$, we get $$(\alpha^{-1}\beta)^*=\beta^*(\alpha^{-1})^*=\alpha^{-1}(\alpha \beta^*)(\alpha^{-1})^*,$$ 
which is $\Lambda$-Hermitian by Lemma \ref{stableunderconjugation}. Hence $\alpha^{-1}\beta$ is $\bar{\Lambda}$-Hermitian. 
Conversely, the condition on $A$ will fulfill the condition $(4)$ of Definition \ref{lamdaunitary}. 
Hence $A$ is $\Lambda$-quadratic. Since $\alpha^{-1}\beta$ is $\bar{\Lambda}$-Hermitian,  
$$T_{12}(-\alpha^{-1} \beta)\in {\EQ}^{\lambda}(2n, R,\Lambda)$$ and $AT_{12}(\alpha^{-1} \beta)= \mbb{H}(\alpha)$. 
Thus $A\equiv \mbb{H}(\alpha) \pmod{{\EQ}^{\lambda}(2n, R,\Lambda)}$. \hb \vp 

A similar proof will prove the following:

\begin{lm}
	\label{lowertriangular}
	Let $B=\begin{pmatrix}
	\alpha &0\\ \gamma &\delta
	\end{pmatrix}\in {\rm M}_{2n}(R)$. Then $B\in {\GQ}^{\lambda}(2n, R,\Lambda)$ if and only if 
	$\alpha\in {\GL}_n(R)$, $\delta=(\alpha^*)^{-1}$ and $\gamma$ is $\Lambda$-Hermitian. 
	In this case  $$B\equiv \mbb{H}(\alpha) \pmod{{\EQ}^{\lambda}(2n, R, \Lambda)}.$$ 
\end{lm}

\begin{lm}
	\label{ifaisinvertible}
	Let $\alpha=\begin{pmatrix}
	a&b\\
	c&d
	\end{pmatrix}\in {\GQ}^{\lambda}(2n,R,\Lambda)$. Then $$\alpha \equiv \mbb{H}(a) \pmod{{\EQ}^{\lambda}(4n,R,\Lambda)}$$ if $a\in {\GL}_n(R).$ Moreover, if  
	$a\in {\E}_n(R)$, 
	then $\alpha\equiv \mbb{H}(a) \pmod{{\EQ}^{\lambda}(2n,R,\Lambda)}$.
\end{lm}

${\pf}$ By same argument as given in Lemma \ref{uppertriangular}, we have $a^{-1}b$ is $\Lambda$-Hermitian. 
Hence  $T_{12}(-a^{-1}b)\in {\EQ}^{\lambda}(2n, R, \Lambda)$, and consequently $\alpha T_{12}(-a^{-1}b)=\begin{pmatrix}
a&0\\c &d^{\prime}
\end{pmatrix}\in {\GQ}^{\lambda}(2n, R, \Lambda)$ for some $d^{\prime}\in {\GL}_n(R)$. Hence by Lemma \ref{lowertriangular}, we get 
$$\alpha T_{12}(-a^{-1}b)\equiv H(a) \pmod{{\EQ}^{\lambda}(2n, R, \Lambda)}.$$ 
Hence  $\alpha\equiv H(a)\pmod{{\EQ}^{\lambda}(2n, R, \Lambda)}$. \hb

\begin{de}\tn{
Let $\alpha=\begin{pmatrix}
a_1&b_1\\c_1&d_1
\end{pmatrix}\in {\rm M}_{2r}(R)$, $\beta=\begin{pmatrix}
a_2&b_2\\c_2&d_2
\end{pmatrix}\in {\rm M}_{2s}(R)$. As before, we define $\alpha\perp \beta$, 
and consider an embedding }
$$\GQ^{\lambda}(2n, R, \Lambda)\rightarrow \GQ^{\lambda}(2n+2,R,\Lambda), \,\,\,\alpha\mapsto \alpha \perp {\I}_2.$$ 
\tn{We denote ${\GQ}^{\lambda}(R,\Lambda)=\underset{n=1}{\overset{\infty}\cup}{\GQ}^{\lambda}(2n, R,\Lambda)$ and 
${\EQ}^{\lambda}(R,\Lambda)=\underset{n=1}{\overset{\infty}\cup} {\EQ}^{\lambda}(2n, R,\Lambda)$.}
\end{de}

In view of quadratic analog of Whitehead Lemma, we have the group ${\EQ}^{\lambda}(R,\Lambda)$ 
coincides with the commutator of ${\GQ}^{\lambda}(R, \Lambda)$. Therefore the group 
$${\k}{\GQ}^{\lambda}(R,\Lambda):= \frac{{\GQ}^{\lambda}(R,\Lambda)}{{\EQ}^{\lambda}(R,\Lambda)}$$ is well-defined. The class of a matrix 
$\alpha\in {\GQ}^{\lambda}(R,\Lambda)$ in the group ${\k}{\GQ}^{\lambda}(R,\Lambda)$ is denoted by $[\alpha]$. In this way we obtain a 
${\k}$-functor ${\k}{\GQ}^{\lambda}$ acting form the category of form rings to the category of abelian groups.

\begin{re}
	\tn{Likewise in the quadratic case, the kernel of the group homomorphism} 
	$${\rm K_1GQ}^{\lambda}(R[X],\Lambda[X])\rightarrow {\rm K_1GQ}^{\lambda}(R,\Lambda)$$
	\tn{induced from  the  form ring homomorphism $(R[X],\Lambda[X])\rightarrow (R,\Lambda); X\mapsto 0$ is denoted by 
		${\rm NK_1GQ}^{\lambda}(R,\Lambda)$. Since the $\Lambda$-quadratic groups are subclass of the quadratic groups, the Local-global principle holds for $\Lambda$-quadratic groups. We use this throughout for the next section.}
\end{re}

\section{\bf Absence of torsion in ${\nk}{\GQ}^{\lambda}(R, \LMD)$}

In this section we give the proof of Theorem \ref{th1} and Theorem \ref{th3}. 
In \cite{RB1}, the proof of the theorem for the linear case is based on 
two key results, {\it viz.} the Higman linearisation, and a lemma on 
polynomial identity in the truncated polynomial rings. 
Here we recall the lemma with its proof to highlight its connection with the 
big Witt vectors. Recently, in \cite{V}, 
V. Kopeiko deduced an analog of Higman linearisation process for a subclass 
of the general quadratic groups.

\begin{de}
	For a associative ring $R$ with unity we consider the truncated polynomial ring 
	$$R_t=\frac{R[X]}{(X^{t+1})}.$$
\end{de}

\begin{lm} \label{la3} $(${\it cf.}\cite{RB1}, Lemma 4.1$)$ 
Let $P(X)\in R[X]$ be any polynomial.
Then the following identity holds in the ring $R_t:$ 
\begin{equation*}
(1+X^r P(X))=(1+X^rP(0))(1+X^{r+1}Q(X)),
\end{equation*}
where $r>0$ and  $Q(X)\in R[X]$, with $\deg(Q(X))< t-r$.
\end{lm}
${\pf}$ Let us write $P(X)=a_0+a_1X+\cdots+a_{t}X^{t}$. Then we can 
write $P(X)=P(0)+XP'(X)$ for some $P'(X)\in R[X]$. Now, in $R_t$
{\small
\begin{align*}
(1+X^r P(X))(1+X^r P(0))^{-1} 
& =  (1+X^r P(0)+X^{r+1}P'(X))(1+X^r P(0))^{-1}\\
& = 1+X^{r+1}P'(X)(1-X^rP(0)+X^{2r}(P(0))^2-\cdots)\\
& =  1+X^{r+1}Q(X)
\end{align*}}
\!\!where $Q(X)\in R[X]$ with $\deg(Q(X))< t-r$. 
Hence the lemma follows.  \hfill$\Box$ \vp \\
{\bf Remark.}  Iterating the above process we can write for any polynomial 
$P(X)\in R[X]$,
$$(1+XP(X))=\Pi_{i=1}^t(1+a_iX^i)$$ in $R_t$, 
for some $a_i\in R$. By ascending induction it will follow that the $a_i$'s 
are uniquely determined. In fact, if $R$ is commutative then 
$a_i$'s are the $i$-th component of the 
ghost vector corresponding to the big Witt vector of 
$(1+XP(X))\in {\w}(R)=(1+XR[[X]])^{\times}$.  For details see 
(\cite{BL1}, $\mc{x}$I).

\begin{lm}
	\label{keylemma}
	Let $R$ be a ring with $1/k\in R$ and $P(X)\in R[X]$. Assume $P(0)$ lies in the center of $R$. Then
	$$(1+X^rP(X))^{k^r}=1 \Rightarrow (1+X^rP(X))=(1+X^{r+1}Q(X))$$ in the ring $R_t$ 
	for some $r>0$ and $Q(X)\in R[X]$ with $\deg(Q(X))<t-r$.
\end{lm}

Following result is due to V. Kopeiko, {\it cf.} \cite{V}.

\begin{pr} {\bf (Higman linearisation)}
\label{kopeiko NK_1} Let $(R, \LMD)$ be a form ring. Then, 
	every element of the group ${\nk}{\GQ}^{\lambda}(R,\Lambda)$ has a 
	representative of the form $$[a;b,c]_n=\begin{pmatrix}
	{\I}_{r}-aX & bX\\-cX^{n}& {\I}_r+a^*X+\cdots +(a^*)^nX^n
	\end{pmatrix}\in {\GQ}^{\lambda}(2r,R[X], \Lambda[X])$$ 
	for some positive integers $r$ and $n$, where $a,b,c\in {\rm M}_r(R)$ satisfy the following conditions:
\begin{enumerate}
 \item the matrices $b$ and $ab$ are Hermitian and also $ab = ba^*$,
 \item the matrices $c$ and $ca$ are Hermitian and also $ca = a^*c$,
 \item  $bc = a^{n+1}$and $cb = (a^*)^{n+1}$.
\end{enumerate}
\end{pr}

\begin{co}
	\label{class is hyperbolic}
	Let $[\alpha]\in {\nk}{\GQ}^{\lambda}(R,\Lambda)$ has the representation  $[a;b,c]_n$ for some $a,b,c\in {\rm M}_n(R)$ according to Proposition  \ref{kopeiko NK_1}. 
	Then $$[\alpha]=[\mbb{H}({\rm I}_r-aX)]$$ in ${\nk}{\GQ}^{\lambda}(R,\Lambda)$ if $({\rm I}_r-aX)\in {\GL}_r(R)$.
	\end{co}

${\pf}$ By Lemma \ref{ifaisinvertible} we have $[a;b,c]_n\equiv \mbb{H}({\I}_r-aX) \pmod{{\EQ}^{\lambda}(2r, R[X],\Lambda[X])}$. 
Hence we have $[\alpha]=[\mbb{H}({\I}_r-aX)]$ in ${\nk}{\GQ}^{\lambda}(R, \Lambda)$. \hb \vp

{\bf Proof of Theorem \ref{th1}:} \vp \\
By the Theorem \ref{kopeiko NK_1}, we have $[\alpha]=[[a;b,c]_n]$ for some $a,b,c\in {\rm M}_s(R)$ and for some natural numbers $n$ and $s$.  
Note that in the Step $1$ of the Proposition \ref{kopeiko NK_1}, the invertibility of the first corner of the matrix $\alpha$ will not be 
changed during the linearisation process.
Also the invertibility of the first corner is preserved in the remaining steps of the Proposition \ref{kopeiko NK_1}. 
Therefore since the first corner matrix $A(X)\in {\GL}_r(R[X])$, then we have $({\I}_s-aX)\in {\GL}_s(R[X])$.
 By Corollary \ref{class is hyperbolic}, we have $[\alpha]=[\mbb{H}({\I}_s-aX)]$. Now let $[\alpha]$ be a $k$-torsion. 
 Then we have $[\mbb{H}({\I}_r-aX)]$ is a $k$-torsion. Since $({\I}_r-aX)$ is invertible, it follows that $a$ is nilpotent. 
 Let $a^{t+1}=0$. Since $[({\I}_r-aX)]^k=[{\I}]$ in ${\k}{\GQ}^{\lambda}(R[X], \Lambda[X])$, 
 then by arguing as given in \cite{RB}, we have $[{\I}_r-aX]=[I]$ in ${\k}{\GQ}^{\lambda}(R[X], \Lambda[X])$. 
 This completes the proof. \hb \vp 

{\bf Proof of Theorem \ref{th3} -- ({\bf Graded Version):} \vp \\
Consider the ring homomorphism $f:R\rightarrow R[X]$ defined by $$f(a_0+a_1+\dots)=a_0+a_1X+\dots.$$ Then
\begin{align*}
[({\I}+N)^k]=[{\I}]&\Rightarrow f([{\I}+N]^k)=[f({\I}+N)]^k=[{\I}]\\
&\Rightarrow [({\I}+N_0+N_1X+\dots +N_rX^r)]^k=[{\I}].
\end{align*}
Let $\mathfrak{m}$ be a maximal ideal in $R_0$.  By Theorem \ref{has no torsion}, 
we have $$[({\I}+N_0+N_1X+\dots+N_rX^r)]=[{\I}]$$ in ${\nk}{\GQ}^{\lambda}((R)_{\mathfrak{m}},\Lambda_{\mathfrak m})$. 
Hence by using the local-global principle we conclude $$[({\I}+N)]=[{\I}+N_0]$$ 
in ${\nk}{\GQ}^{\lambda}(R,\Lambda)$, as required. \hb \vp\\
{\bf Acknowledgment:} We thank Sergey Sinchuk and V. Kopeiko for many useful discussions.

{\small

\medskip

\end{document}